\documentclass{amsart}
 
\usepackage{latexsym, amssymb, amsfonts, amsmath, amsthm, graphics, graphicx, subfigure, bbm, stmaryrd}
\usepackage{epstopdf}

\newcommand{\field}[1]{\mathbb{#1}}
\newcommand{\R}{\field{R}}
\newcommand{\E}{\field{E}}
\newcommand{\Z}{\field{Z}}
\newcommand{\C}{\field{C}}

\newtheorem{theorem}[equation]{Theorem}

\newtheorem{proposition}[equation]{Proposition}

\newtheorem{lemma}[equation]{Lemma}
\newtheorem{corollary}[equation]{Corollary}
\newtheorem{remark}[equation]{Remark}

\newcommand{\met}[2][ccccccccccccccccccccccccccccccccc]{\left[ \begin{array}{#1} #2 \\ \end{array}\right]}

\numberwithin{equation}{section}

\title[ Talagrand Inequality for the Semicircular]{\bf{Talagrand Inequality for the Semicircular Law and Energy of the Eigenvalues of Beta Ensembles}}
\author{Ionel Popescu}
\date{}

\address{Department of Mathematics\\  Northwestern University\\  2033 Sheridan Road, Evanston, IL, 60208, USA}
\address{IMAR, 21 Calea Grivitei Street, 010702-Bucharest, Sector 1, Romania} 

  \email{ipopescu@math.northestern.edu}

\begin{document}

\begin{abstract} 
We give a short proof of an extension of the free  Talagrand  transportation cost inequality to the semicircular which was originally proved in \cite{BV}.   The proof is based on a convexity argument and is in the spirit of the original Talagrand's approach for the classical counterpart from \cite{T}.  We also discuss the convergence, fluctuations and large deviations of the energy of the eigenvalues of $\beta$ ensembles, which, as an application of Talagrand inequality gives in particular yet another proof of the convergence of the eigenvalue distribution to the semicircle law.
\end{abstract}

\maketitle

\section{Introduction}

In \cite{T}, Talagrand proves the transportation cost inequality to the Gaussian measure. The one dimensional version for the Gaussian measure  $\gamma(dx)=\frac{1}{\sqrt{2\pi}}e^{-x^{2}/2}dx$ reads as
\begin{equation}\label{0}
(W_{2}(\mu,\gamma))^{2}\le 2 H(\nu |\gamma), 
\end{equation}
where $W_{2}(\mu,\gamma)$ is the Wasserstein distance defined below by \eqref{e4} and the relative entropy  is
\[
H(\nu | \gamma) = 
\begin{cases}
\int f(x)\log(f(x))d\gamma(x) & \text{if}\: \nu(dx)=f(x)\gamma(dx) \\  
\infty &\text{if}\: \nu\:\text{is singular to}\: \gamma.  
\end{cases}
\]

 In the context of free probability, Biane and Voiculescu proved in \cite{BV} a free version of this: 
 \begin{equation}\label{1}
 (W_{2}(\mu,\sigma))^{2}\le  2(E(\mu) - E(\sigma)),
 \end{equation}
where $E(\mu)=\frac{1}{2} \int x^{2}\mu(dx) -\iint \log(|x-y|)\mu(dx)\mu(dy)$  is the free energy of $\mu$ and $\sigma(dx)=\frac{1}{2\pi}\mathbbm{1}_{[-2,2]}(x)\sqrt{4-x^{2}}dx$ is the semicircular law, the minimizer of $E(\mu)$ over all probability measures on the real line.   The role of the relative entropy is played here by the difference of the free energy of $\mu$ and the semicircular.  

Using random matrix approximations, Hiai, Petz and Ueda  proved in \cite{HPU} the following extension of \eqref{1},
\begin{equation}\label{2}
   \rho \left(W_{2}(\mu, \mu_{Q})\right)^{2} \le E^{Q}(\mu) -E^{Q}(\mu_{Q})
\end{equation}
where $\rho>0$ and $Q:\R\to\R$ is a function so that $Q(x)-\rho x^{2}$ is convex and 
\[
E^{Q}(\mu)=\int Q(x)\mu(dx)-\iint\log|x-y|\mu(dx)\mu(dy).
\]  
Here $\mu_{Q}$ is the minimizer of $E^{Q}$ on the set of all probability measures on the real line. They also prove a version of this for measures supported on the circle $\mathbf{T}$: 
\begin{equation}\label{3}
   (\rho+1/4) \left(W_{2}(\nu, \nu_{Q})\right)^{2} \le E^{Q}(\nu) -E^{Q}(\nu_{Q})
\end{equation}
where $Q:\mathbf{T}\to\R$ so that $Q(e^{i x})-\rho x^{2}$ is convex on $\R$, $\rho>-1/4$ and $\mu_{Q}$ is the minimizer of the functional $E^{Q}$ on probability measures on the unit circle $\mathbf{T}$.  

Another proof of \eqref{1} is given in \cite{L} via a Brunn-Minkovsky inequality for free probability.  

The primary purpose of this note is to give an elementary proof of \eqref{2}  and \eqref{3} in the spirit of  Talagrand's proof to \eqref{0}.  The idea is to exploit convexity of the logarithm appearing in the $E^{Q}$. We also discuss (see Theorem \ref{TQd} and Proposition \ref{prop}) the discrete version of the transportation cost inequalities and some consequences involving Fekete points.  

The second purpose of this note is to discuss the energy of the eigenvalues of $\beta$ ensembles and in particular the fluctuations and the deviations from the minimum energy (see Theorem~\ref{B}).   This is a simple application of Selberg's formula together with elementary estimates on $\Gamma$ functions.  As a consequence, using the the results in the first part we reprove that  the distribution of the eigenvalues converges almost surely to the semicircular law.

\section{Talagrand Inequalities}

The following result is an obvious one but is the key to our problem.
\begin{lemma}\label{l1} Let $f:[0,1]\to\R$ be a convex function with the property that $f(0)=0$ and there exists $a\ge0$ so that
\[
f(t) \ge -a t^{2}\quad \text{for}\quad t\in [0,1].
\] 
Then
\[
f(t)\ge0\quad \text{for all}\quad t\in [0,1].
\]
\end{lemma}

\begin{proof} 
It follows from the assumptions that for any $\epsilon>0$, if $\delta_{\epsilon}=\min(1,\epsilon/a)$, then $f(t)\ge-t\epsilon $ for $t\in[0,\delta_{\epsilon}]$.  Now, since $f$ is convex, one gets $f(mt)\ge mf(t)\ge-mt\epsilon$ for any integer $m$ with $mt\le 1$, and  therefore, $f(t)\ge - \epsilon t$ for any $t\in[0,1]$.  Since this is true for any $\epsilon > 0$, we get $f(t)\ge 0$ for any $t\in[0,1]$.
 \end{proof}
 
 In the following,  $\mathcal{P}(\Omega)$ denotes the set of all probability measures on $\Omega$, and for two probability measures  with finite second moment on $\mathcal{P}(\R)$ or $\mathcal{P}(\mathbf{T})$, where $\mathbf{T}=\{ z\in\C:|z|=1\}$, we define $W_{2}(\mu,\nu)$,  the Wasserstein distance by
 \begin{equation}\label{e4}
 W_{2}(\mu,\nu):=\sqrt{\inf_{\pi\in\Pi(\mu,\nu)} \iint|x-y|^{2}d\pi(x,y)}.
 \end{equation}
 Here $\Pi(\mu,\nu)$ is the set of probability measures on $\R^{2}$ with marginal distributions $\mu$ and $\nu$, and  it can be shown that there is at least one solution $\pi\in\Pi(\mu,\nu)$ to this minimization problem.  
 
 If $\mu$ and $\nu$ are two measures on $\R$ with $F$ and $G$ their cumulative distribution functions (i.e. $F(x)=\mu((-\infty,x])$), then Theorem 2.18 in \cite{Vil} states that
 \begin{equation}\label{Wf}
  (W_{2}(\mu,\nu))^{2}=\int_{0}^{1}|F^{-1}(t)-G^{-1}(t)|^{2}dt
  \end{equation}
  where $F^{-1}$ denotes the generalized inverse of $F$.  
  
 \begin{theorem}\label{TQ} Let $Q:\R\to\R$ be a function so that $Q(x)-\rho x^{2}$ is
 convex for a certain $\rho>0$.  If $\mu_{Q}$ is a solution to the minimization problem
 \begin{equation}\label{e2}
I^{Q} := \inf_{\mu\in\mathcal{P}(\R)} E^{Q}(\mu),
 \end{equation}
 where 
 \begin{equation}\label{E}
 E^{Q}(\mu)=\int Q(x)\mu(dx)-\iint\log|x-y|\mu(dx)\mu(dy),
 \end{equation}
  then for any $\mu\in\mathcal{P}(\R)$, we have 
\begin{equation}\label{e5}
  \rho \left(W_{2}(\mu, \mu_{Q})\right)^{2} \le E^{Q}(\mu) -I^{Q}.
\end{equation}
In particular, the minimization problem \eqref{e2} has a unique solution.
 \end{theorem}

 \begin{proof}
 There exist constants $c_{1}$ and $c_{2}$ so that 
 \[
 Q(x)-\rho x^{2} \ge c_{1}\quad \text{and} -\log(|x-y|)\ge -\frac{\rho}{4}(x^{2}+y^{2})+c_{2}.
 \] 
 Then for a certain $C$, we get that 
 \begin{equation}\label{e3}
 \frac{1}{2}(Q(x)+Q(y))-\log(|x-y|) \ge \frac{\rho}{4}(x^{2}+y^{2})+C \ge C,
 \end{equation} 
 and this in turn implies that the infimum in \eqref{e2} is finite (since $E^{Q}(\mu)$ is finite for $\mu$ the uniform distribution on $[0,1]$) and in particular $\int Q(x)d\mu_{Q}(x)$, and $\iint \log|x-y|d\mu_{Q}(x)d\mu_{Q}(y)$ are finite, which means that $\mu_{Q}$ has finite second moment and no atoms.  
 
Since $E^{Q}(\mu)>-\infty$, we may assume that $E^{Q}(\mu)$ is finite, otherwise there is nothing to prove.  Then, $\iint\log|x-y|\mu(dx)\mu(dy)$ and $\int Q(x)\mu(dx)$ are finite.  
In particular, $\mu$ has finite second moment and no atoms.  

Taking $F_{\mu}$ and $F_{\mu_{Q}}$, the cumulative distributions of $\mu$, $\mu_{Q}$ and $F^{-1}$, $F^{-1}_{Q}$ their generalized inverses, set $\theta(x)=F^{-1}(F_{Q}(x))$.   According to \cite[Theorem 2.18]{Vil} and the discussion following thereafter, the minimizing measure $\pi$ from \eqref{e4} is the distribution of $x\to(x,\theta(x))$ under $\mu_{Q}$.   In this case, the inequality we want to prove becomes
\[
  \rho \iint |x-\theta(x)|^{2}\mu_{Q}(dx) \le \int Q(x)\mu(dx)-\iint\log|x-y|\mu(dx)\mu(dy)-I^{Q}.
\]

Let $f:[0,1]\to\R$ be given by
\[
\begin{split}
f(t)=-&\rho t^{2}\int |\theta(x)-x|^{2}\mu_{Q}(dx) + \int Q(t\theta(x)+(1-t)x)\mu_{Q}(dx)\\ 
&-\iint\log(|t(\theta(x)-\theta(y))+(1-t)(x-y)|)\mu_{Q}(dx)\mu_{Q}(dy) -I^{Q}.
\end{split}
\]
Notice here that $f$ is well defined.  Indeed, $Q$ is convex, hence bounded below and  because $\int Q(\theta(x))\mu_{Q}(dx) = \int Q(x)\mu(dx)$ and $\int Q(x)\mu_{Q}(dx)$ are both finite, one concludes that $ \int Q(t\theta(x)+(1-t)x)\mu_{Q}(dx)$ is finite too.  One the other hand, there is a  $C>0$ so that for any $t\in[0,1]$, 
\[
 -\log(|t(\theta(x)-\theta(y))+(1-t)(x-y)|)\ge -C(\theta(x)^{2}+\theta(y)^{2}+x^{2}+y^{2}) - C,
 \]
 which, combined with the finiteness of the second moment of $\mu$ and $\mu_{Q}$, results with (for a constant $C$)
 \[
 -\iint\log(|t(\theta(x)-\theta(y))+(1-t)(x-y)|)\mu_{Q}(dx)\mu_{Q}(dy)>C\quad\text{for all}\quad t\in[0,1].
 \] 
Now, since $\theta$ is a nondecreasing function we can write
\[
\begin{split}
-\iint\log(&|t(\theta(x)-\theta(y))+(1-t)(x-y)|)\mu_{Q}(dx)\mu_{Q}(dy)=\\
 & -2\iint_{x>y}\log(t(\theta(x)-\theta(y))+(1-t)(x-y))\mu_{Q}(dx)\mu_{Q}(dy),
\end{split}
\]
which combined with the convexity of $-\log$ on $(0,\infty)$ and the finiteness of $\iint \log|x-y|\mu_{Q}(dx)\mu_{Q}(dy)$ and $\iint\log|x-y|\mu(dx)\mu(dy)$, yields the fact that 
\[\tag{**}
t\to  -\iint\log(|t(\theta(x)-\theta(y))+(1-t)(x-y)|)\mu_{Q}(dx)\mu_{Q}(dy)
\]
is well defined and convex.

The inequality \eqref{e5} is now equivalent to $f(1)\ge0$.  To show this, we apply Lemma~\ref{l1}.  
 The convexity follows easily from the convexity of $Q(x)-\rho x^{2}$ and (**).   Now if $\nu_{t}$ is the distribution of $x\to t\theta(x)+(1-t)x$ under $\mu_{Q}$, then the minimization property of $\mu_{Q}$ implies that
\[
f(t) \ge -\rho t^{2}\iint |\theta(x)-x|^{2}\mu_{Q}(dx)\quad\text{for}\quad t\in[0,1],
\]
and then, Lemma~\ref{l1} shows that $f(t)\ge 0$ for any $t\in[0,1]$.

The existence statement follows from the lower continuity of $E^{Q}$.  For a proof of the existence and compactness of the support of $\mu_{Q}$, see for instance Chapter $6$ in \cite{D}. 
 \end{proof}

 \begin{remark}  What was essential during the proof was the convexity of $-\log$ on $(0,\infty)$ and the fact that for any $a>0$, there is a $C(a)$ so that $-\log|x-y|\ge -a(x^{2}+y^{2})+C(a)$.  Therefore if we replace the $\log$ in the statement of this theorem by any kernel $K(|x-y|)$ with the property that 
 $K$ on $(0,\infty)$ is concave and that for any $a>0$, there is a $C(a)$ so that $-K(|x-y|)\ge -a(x^{2}+y^{2})+C(a)$, then the result still holds.   Other examples of such kernels are $-1/x^{\alpha}$, $\alpha>0$ and $1/\log(x^{2}+1)$.
 \end{remark}

 If we take $Q(x)=\frac{x^{2}}{2}$, and keep in mind that the minimizing measure $\mu_{Q}$ for $E^{Q}$ is the semicircular law, one gets the following result proved in \cite{BV}.
 \begin{corollary}\label{cor}
 Let $\sigma(dx)=\frac{1}{2\pi}\mathbbm{1}_{[-2,2]}(x)\sqrt{4-x^{2}}dx$ be the semicircular law on $[-2,2]$.  Then for any $\mu\in\mathcal{P}(\R)$,
 \[
\frac{1}{2} (W_{2}(\mu,\sigma))^{2}\le \frac{1}{2} \int x^{2}\mu(dx) -\iint \log(|x-y|)\mu(dx)\mu(dy) -\frac{3}{4}.
 \]
 \end{corollary}
 
 The next theorem is just inequality \eqref{3}.
 
 \begin{theorem}
 Assume $Q:\mathbf{T}\to\R$ is a function  so that $Q(e^{ix})-\rho x^{2}$ is convex on $\R$ for a given $\rho>-1/4$.  If $\mu_{Q}$ is a solution to the minimization problem 
 \begin{equation}\label{e2'}
I^{Q} := \inf_{\mu\in\mathcal{P}(\mathbf{T})} E^{Q}(\mu),
 \end{equation}
 where 
 \begin{equation}\label{E'}
 E^{Q}(\nu)=\int Q(z)\nu(dz)-\iint_{\mathbf{T}\times\mathbf{T}}\log|z-z'|\nu(dz)\nu(dz'),
 \end{equation}
 then, for any $\nu\in\mathbf{T}$, we have 
\begin{equation}\label{e5'}
(\rho+1/4) \left(W_{2}(\nu, \nu_{Q})\right)^{2} \le E^{Q}(\nu) -I^{Q}.
\end{equation} 
In particular, there is a unique solution for the minimization problem \eqref{e2'}.
 \end{theorem}
 
 \begin{proof}
 Take the exponential map $\exp:x\in\R\to e^{ix}\in\mathbf{T}$ and for any measure $\mu$ on $\mathbf{T}$, define $\bar{\mu}(A)=\sum_{n\in\Z}\mu(\exp(A\cap [2\pi n,2\pi(n+1))))$.  One can show that there exists $L\in[0,2\pi)$ such that the restrictions of $\bar{\mu}$ and $\bar{\mu}_{Q}$ to $[L,L+2\pi)$ have the same mean value.  
  
 We then identify $[L,L+2\pi)$ with $\mathbf{T}$ via the exponential map and define $\nu_{Q}$ and $\nu$ to be the restrictions of $\bar{\mu}_{Q}$ and $\bar{\mu}$ to the interval $[L,L+2\pi)$.    We then follow the proof of \ref{TQ} with the necessary adjustments needed.  We take the function $f(t)$ here to be 
 \[
\begin{split}
f(t)=-&(\rho+1/4) t^{2}\int |\theta(x)-x|^{2}\nu_{Q}(dx) + \int Q(e^{i(t(\theta(x)+c)+(1-t)x)})\nu_{Q}(dx)\\ 
&-\iint\log(|e^{i(t(\theta(x)+(1-t)x))}-e^{i(t\theta(y)+(1-t)(y))}|)\nu_{Q}(dx)\nu_{Q}(dy) -I^{Q}.
\end{split} 
 \]
Now, $|e^{ia}-e^{ib}|^{2}=4\sin^{2}( (a-b)/2)$ for  $a$, $b$ real numbers and 
\[
\int |\theta(x)-x|^{2}\nu_{Q}(dx)=\frac{1}{2}\iint((\theta(x)-x)-(\theta(y)-y))^{2}\nu_{Q}(dx)\nu_{Q}(dy).
\]
Next, set $\theta_{t}(x)=t\theta(x)+(1-t)x$ and notice that
\begin{align*}
g(t):=&- \frac{t^{2}}{4}  \int |\theta(x)-x|^{2}\nu_{Q}(dx)- \iint\log(|e^{it\theta_{t}(x))}-e^{i\theta_{t}(y))}|)\nu_{Q}(dx)\nu_{Q}(dy) \\
 = & - \iint\frac{t^{2}}{8}((\theta(x)-x)-(\theta(y)-y))^{2}\nu_{Q}(dx)\nu_{Q}(dy) \\ & \qquad\qquad-\iint \log\left|2\sin\left((\theta_{t}(x)-\theta_{t}(y))/2\right)\right|\nu_{Q}(dx)\nu_{Q}(dy)\\ 
 =  & - 2\iint_{x>y}\frac{t^{2}}{8}((\theta(x)-x)-(\theta(y)-y))^{2}\nu_{Q}(dx)\nu_{Q}(dy) \\ 
 &\qquad \qquad-2\iint_{x>y} \log \left(2\sin\left((\theta_{t}(x)-\theta_{t}(y))/2\right)\right)\nu_{Q}(dx)\nu_{Q}(dy),
\end{align*}
where in the last line we used the fact that $\theta$ is a nondecreasing function.  Since $x,y,\theta(x),\theta(y)\in[-\pi,\pi)$ and for $0<a<b<\pi$, we have
\begin{align*}
\frac{d^{2}}{dt^{2}}&\left( -\frac{t^{2}}{8}(a-b)^{2}-\log\left(\sin\left(\frac{ta + (1-t) b}{2}\right)\right) \right)\\ & =\frac{(a-b)^{2}}{4}\left( \frac{1}{\sin^{2}\left(\frac{ta + (1-t) b}{2}\right)}-1 \right)\ge 0,
\end{align*}
which implies that the function $g$ is convex on $[0,1]$.  This coupled with the convexity of $Q(e^{ix})-\rho x^{2}$ concludes that $f$ is a convex function.  Finally
\[
f(t)\ge - (\rho+1/4)t^{2}\int |\theta(x)-x|^{2}\nu_{Q}(dx),
\]
and thus, Lemma~\ref{l1} shows that $f(1)\ge0$, which is \eqref{e5'}.

The existence of a minimizer follows from the fact that $E^{Q}$ is lower semicontinuous.  
\qedhere
\end{proof}

For $Q=0$ and $\rho=0$, the minimizer of \eqref{e2'} is the Haar measure on $\mathbf{T}$.  One can check this by showing directly that the uniform measure satisfy the variational form of \eqref{e2'}. 

\begin{corollary}
For any $\mu\in\mathcal{P}(\mathbf{T})$
\[
\frac{1}{4}\left(W_{2}\left(\mu,\frac{dx}{2\pi}\right)\right)^{2}\le -\iint_{\mathbf{T}\times\mathbf{T}}\log|z-z'|\mu(dz)\mu(dz').
\]
\end{corollary}

 Using the same argument as in the proof of Theorem~\ref{TQ}, we can also prove a discrete version of it. 
 
 \begin{theorem}\label{TQd}
 Let $Q:\R\to\R$ be a function so that $Q(x)-\rho x^{2}$ is 
 convex for a certain $\rho>0$.  For $\mathbf{x}=(x_{1},x_{2},\dots,x_{n})\in\R^{n}$, set the energy of $\mathbf{x}$ to be given by
 \[
 E^{Q}_{n}(\mathbf{x})=\frac{1}{n}\sum_{k=1}^{n}Q(x_{i})-\frac{2}{n(n-1)}\sum_{1\le i<j\le n}\log|x_{i}-x_{j}|.
 \]
If $\Delta^{Q}_{n}=E^{Q}_{n}(\mathbf{y})=\inf\{ E^{Q}_{n}(\mathbf{x}):x\in\R^{n} \}$, then for any $\mathbf{x}\in\R^{n}$, 
\begin{equation}\label{dT}
\rho(W_{2}(\mu(\mathbf{x}),\mu(\mathbf{y})))^{2}\le E^{Q}_{n}(\mathbf{x})-E^{Q}_{n}(\mathbf{y})=E^{Q}_{n}(\mathbf{x})-\Delta^{Q}_{n}
\end{equation}
where $\mu(\mathbf{x})=\frac{1}{n}\sum_{k=1}^{n}\delta_{x_{k}}$.  Moreover, 
\begin{equation}\label{inc}
\Delta^{Q}_{n}\le \Delta_{n+1}^{Q}.
\end{equation}
\end{theorem}

The only statement that needs to be clarified here is \eqref{inc}.  If $\mathbf{y}_{n+1}$ is a minimum point for $E^{Q}_{n+1}$ and  $\mathbf{y}_{n+1}^{i}$ denotes the $n$ dimensional vector obtained from $\mathbf{y}_{n+1}$ by removing the  $i$th component, then $\Delta_{n+1}^{Q}=\frac{1}{n+1}\sum_{i=1}^{n+1}E^{Q}_{n}(\mathbf{y}_{n+1}^{i})$, which is obviously $\ge \Delta_{n}^{Q}$.  
 
The minimum points of $E^{Q}_{n}$ are called Fekete points in the literature.  It is known (see for instance chapter $6$ in \cite{D}) that $\lim_{n\to\infty}\Delta_{n}^{Q}=I^{Q}$, with $I^{Q}$ defined in \eqref{e2}.   We will reprove this fact below in Proposition~\ref{prop}.

For $Q(x)=x^{2}$, the formula \cite[A.6.11]{M} with the appropriate scaling gives the formula for computing $\Delta_{n}=\Delta_{n}^{Q}$ as
\begin{equation}\label{e8}
\Delta_{n}=\frac{1}{2}(1+\log(n-1))-\frac{1}{n(n-1)}\sum_{j=1}^{n}j\log j = \frac{1}{2}-\frac{\log n}{n-1}-\frac{1}{n}\sum_{j=1}^{n-1}\frac{j}{n-1}\log\left( \frac{j}{n-1}\right). 
\end{equation}
 
The next statement is a similar result to Theorems \ref{TQ} and \ref{TQd}.
 
 \begin{proposition}\label{prop}
 Assume $Q:\R\to\R$ is a function so that $Q(x)-\rho x^{2}$ is  convex for a certain $\rho>0$.  Then for any $\nu\in\mathcal{P}(\R)$ and $\mathbf{y}\in\R^{n}$ a Fekete point for 
 $E^{Q}_{n}$, we have 
 \begin{equation}\label{e11}
 \rho (W_{2}(\nu,\mu(\mathbf{y})))^{2}\le E^{Q}(\nu)-\Delta^{Q}_{n}.
 \end{equation}
 Furthermore, if $\mu_{Q}$ is the minimizing measure of $E^{Q}$, and $\mathbf{y}_{n}\in\R^{n}$ is a Fekete point for $E^{Q}_{n}$, then
 \begin{equation}\label{e10} 
 \lim_{n\to\infty}\Delta^{Q}_{n}=I^{Q}\quad\text{and}\quad \lim_{n\to\infty}W_{2}(\mu_{Q},\mu(\mathbf{y}_{n}))=0,
 \end{equation}
 hence, $\mu(\mathbf{y}_{n})\underset{n\to\infty}{\longrightarrow}\mu_{Q}$ weakly.  
 \end{proposition}

\begin{proof} In the first place there is nothing to prove if $E^{Q}(\nu)=\infty$.  Therefore we assume that $E^{Q}(\nu)<\infty$.  Integrating \eqref{dT} with respect to $\nu(dx_{1})\nu(dx_{2})\dots\nu(dx_{n})$, one gets that 
\[
\rho\int (W_{2}(\mu(\mathbf{x}),\mu(\mathbf{y})))^{2}\nu(dx_{1})\nu(dx_{2})\dots\nu(dx_{n}) \le E^{Q}(\nu)-\Delta^{Q}_{n}. 
\]
 We finish the proof of \eqref{e11} by showing that
\[\tag{*}
\int (W_{2}(\mu(\mathbf{x}),\mu(\mathbf{y})))^{2}\nu(dx_{1})\nu(dx_{2})\dots\nu(dx_{n})=(W_{2}(\nu,\mu(\mathbf{y})))^{2}.
\]
To do this, we proceed by induction.  For $n=1$, this statement becomes
\[
\int (W_{2}(\delta_{x},\delta_{y}))^{2}\nu(dx)=(W_{2}(\nu,\delta_{y}))^{2}
\]
which, cf. \eqref{Wf}, is equivalent to the following (here $F_{\nu}$ is the cumulative distribution function of $\nu$)
\[
\int |x-y|^{2}\nu(dx)=\int_{0}^{1}|y-F_{\nu}^{-1}(t)|^{2}dt.
\]
This can be checked by changing the variable in the second integral.

Assume (*) is true for $n-1$, $n\ge2$.  A simple application of \eqref{Wf} gives that 
$(W_{2}(\mu(\mathbf{x}),\mu(\mathbf{y})))^{2}=\frac{1}{n}\sum_{i=1}^{n}|x_{\sigma(i)}-y_{\tau(i)}|^{2}$, where $\sigma$ and $\tau$ are permutations of $\{1,2,\dots,n \}$ so that $x_{\sigma(1)}\le x_{\sigma(2)}\dots \le x_{\sigma(n)}$ and $y_{\tau(1)}\le y_{\tau(2)}\dots \le y_{\tau(n)}$.  If we denote by $\mathbf{x}_{i}$ the vector $\mathbf{x}$ with the $i$th component removed and similarly for $\mathbf{y}_{i}$, one deduces 
\[\tag{\#}
(W_{2}(\mu(\mathbf{x}),\mu(\mathbf{y})))^{2}=\frac{1}{n}\sum_{i=1}^{n}(W_{2}(\mu(\mathbf{x}_{i}),\mu(\mathbf{y}_{i})))^{2}.
\]  
On the other hand,
\[
(W_{2}(\nu,\mu(\mathbf{y})))^{2}=\sum_{k=0}^{n-1}\int_{k/n}^{(k+1)/n}|y_{\tau(k)}-F_{\nu}^{-1}(t)|^{2}dt,
\]
which can be used to argue that
\[\tag{\#\#} 
(W_{2}(\nu,\mu(\mathbf{y})))^{2}= \frac{1}{n}\sum_{i=1}^{n}(W_{2}(\nu,\mu(\mathbf{y}_{i})))^{2}.
\]   
 Putting together ($\#$) and ($\#\#$) and the induction hypothesis one finishes the proof of (*).    

To prove \eqref{e10}, we first point out that \eqref{e11} applied to $\mu_{Q}$ yields that $I^{Q}\ge \Delta_{n}^{Q}$ for any $n\ge1$.  In particular this means that $\Delta_{n}^{Q}$ is bounded.  Since $-\log|x-y|\ge -\frac{\rho}{4}(x^{2}+y^{2})+c$ for a certain constant $c$, we get that $\Delta_{n}^{Q}\ge \frac{\rho}{4n}\sum_{i=1}^{n}x_{i}^{2}-C$, where $C$ is a constant.  This implies that the sequence $\{\int x^{2}\mu(\mathbf{y}_{n})(dx)\}_{n\ge1}$ is bounded, whose consequence is that the sequence of measures $\mu(\mathbf{y}_{n})$ is tight, therefore there is a weak convergent subsequence $\mu(\mathbf{y}_{n_{k}})$ to a measure $\nu$.  Now, for any $L>0$, we have 
\[
\int \min\{ ((Q(x)+Q(y))/2-\log|x-y|),L\}\mu(\mathbf{y}_{n_{k}})(dx)\mu(\mathbf{y}_{n_{k}})(dy)\le \Delta_{n_{k}}^{Q}+L/n_{k}
\]
and this demonstrates that for any $L>0$,
\[
\int \min\{ ((Q(x)+Q(y))/2-\log|x-y|),L\}\nu(dx)\nu(dy)\le I^{Q}
\]
and, after passing $L\to\infty$, this yields 
\[
E^{Q}(\nu)\le I^{Q}.
\]
This together with \eqref{inc} and the uniqueness of $\mu_{Q}$ from Theorem~\ref{TQ}  ends the proof of $\lim_{n\to\infty}\Delta_{n}^{Q}=I^{Q}$.  The rest follows.  \qedhere
\end{proof}

 \section{Discrete Energy for $\beta$-Ensembles}
 
 In this section we deal with $\beta$-ensembles, which are studied in \cite{DE}.   These are tridiagonal matrices with independent entries of the form
 \[
 A_{n}=\frac{1}{\sqrt{\beta n}} 
 \met{ N(0,2) & \chi_{(n-1)\beta} & & & \\ \chi_{(n-1)\beta} & N(0,2) & \chi_{(n-2)\beta} & & \\ 
 & \ddots & \ddots & \ddots & \\ & & \chi_{2\beta} & N(0,2) & \chi_{\beta} \\ & & &\chi_{\beta} & N(0,2) }.
 \]
Here $N(0,2)$ stands for a normal with mean $0$ and variance $2$, while $\chi_{\gamma}$ is the $\chi$-distribution with parameter $\gamma$.  The joint distribution of the eigenvalues is 
\[
\frac{1}{Z_{\beta,n}}\prod_{1\le i<j\le n}|x_{i}-x_{j}|^{\beta}\exp\left(- \beta n\sum_{i=1}^{n}x_{i}^{2} \right)
\]
where here $Z_{\beta,n}$ is a normalization constant.  

Set $\mu_{n}=\sum_{k=1}^{n}\delta_{\lambda_{k,n}}$, the empirical distribution of the eigenvalues $\{\lambda_{k,n}\}_{k=1}^{n}$ of $A_{n}$.  

\begin{theorem}\label{B}
Set $E_{n}=\frac{1}{2n}\sum_{k=1}^{n}\lambda_{k}^{2}-\frac{2}{(n-1)n}\sum_{1\le j<k\le n}\log|\lambda_{i}-\lambda_{j}|$ the energy of the eigenvalues $\{\lambda_{k}\}_{k=1}^{n}$ of $A_{n}$.  If $\Delta_{n}$ is the quantity defined in \eqref{e8}, then almost surely,
\begin{equation}\label{con}
\lim_{n\to\infty}n(E_{n}-\Delta_{n}) = \psi(1+\beta/2)-\log(\beta/2),
\end{equation}
where $\psi(x)=\frac{d}{dx}\log \Gamma(x)$ and $\Gamma$ is the Gamma function.  In addition, we have that
\begin{equation}\label{dis}
n^{1/2}\bigg(n(E_{n}-\Delta_{n})-(\psi(1+\beta/2)-\log(\beta/2))\bigg) \xrightarrow[n\to\infty]{}  N(0,\psi'(1+\beta/2)),
\end{equation}
where the convergence is in distribution sense.  

The large deviations of $n(E_{n}-\Delta_{n})$ is governed by the rate function
\[
R^{*}(t)=\sup\{tz-R(z):z\in\R \},
\]
\[
R(z)=
\begin{cases}
z+(\beta/2-z)\log(\beta/2-z)-\log\left(\frac{\Gamma(1+\beta/2-z)}{\Gamma(1+\beta/2)}\right)-(\beta/2)\log(\beta/2), & z<\beta/2 \\ 
\infty, & z\ge\beta/2. 
\end{cases}
\]
\end{theorem}

\begin{proof}
The proof is based on a version of Selberg's formula and elementary approximations involving Gamma function.  

First, we have
\[
\E\left[\exp(z \E_{n}) \right] = \frac{\int_{\R^{n}} \prod_{1\le i<j\le n}|x_{i}-x_{j}|^{\beta-\frac{2z}{n(n-1)}}\exp\left(-(\beta n -\frac{z}{2n})\sum_{j=1}^{n}x_{j}^{2} \right)d\mathbf{x} }{\int_{\R^{n}} \prod_{1\le i<j\le n}|x_{i}-x_{j}|^{\beta}\exp\left(-\beta n\sum_{j=1}^{n}x_{j}^{2} \right)d\mathbf{x}} 
\]     
and then, as a consequence of Selberg's formula \cite[equation 17.6.7]{M}, we get for complex $z$, that
\[
\E\left[e^{z \E_{n}} \right]=
\begin{cases}
\frac{\left( n\beta/2 -z/n\right)^{-\frac{n}{2}\left[(n-1)(\beta/2-\frac{z}{n(n-1)}+1)\right]} \prod_{j=1}^{n}\frac{\Gamma(1+j(\beta/2-\frac{z}{n(n-1)}))}{\Gamma(1+\beta/2)}}{(n\beta/2)^{-\frac{n}{2}[(n-1)\beta/2+1]}\prod_{j=1}^{n}\frac{\Gamma(1+j\beta/2)}{\Gamma(1+\beta/2)}}, &\Re(z)< \beta/2 \\ 
\infty, &\Re(z)\ge\beta/2.
\end{cases}
\]
We need Stirling formula for approximation of Gamma function in the following form 
\[
\log\Gamma(t+1)=(t+1/2)\log t - t +\log (2\pi)/2+\mathcal{O}\left(\frac{1}{1+t}\right)\:\text{for}\:t\ge0
\]
Using this and the above formula for $\E[\exp(zE_{n})]$ and \eqref{e8}, after some arrangements one gets 
\begin{equation}\label{k}
\begin{split}
&\log(\E\left[e^{z(E_{n}-\Delta_{n})}\right]) =\frac{z}{n-1}+\frac{z}{2}\log\left(1+\frac{1}{n-1} \right)\\  & -\frac{z}{n-1}\log\left( \frac{\beta}{2} -\frac{z}{n^{2}}\right)+\frac{z(n+1)}{2(n-1)}\log\left( 1+\frac{z}{n[(n-1)n\beta/2-z]} \right)  \\ &+  \frac{n[(n-1)\beta+1]}{2}\log\left(1-\frac{z}{(n-1)[n^{2}\beta /2-z]} \right)+\frac{n\beta}{2}\log\left( 1-\frac{2z}{n(n-1)\beta} \right)\\ & -n\left[ \log\left( 1+\frac{\beta}{2}-\frac{z}{n(n-1)} \right)-\log\left( 1+\frac{\beta}{2}\right)\right] +\mathcal{O}\left(\frac{z}{n^{2}}\right).
\end{split}
\end{equation}    
From this, replacing $z$ by $nz$, one immediately obtains that for any $z\in\R$,   
\[
\log(\E[\exp(zn(E_{n}-\Delta_{n}))]) \underset{n\to\infty}{\longrightarrow} z\frac{\Gamma'(1+\beta/2)}{\Gamma(1+\beta/2)}-z\log(\beta/2)=z(\psi(1+\beta/2)-\log(\beta/2)).
\]
Applying\eqref{k} with $z$ replaced by $n^{3/2}z$, one can prove that for any complex $z$,  
\[
\log\left(\E\left[ \exp\left( zn^{1/2}(n(E_{n}-\Delta_{n})-(\psi'(1+\beta/2)-\log(\beta/2)) \right) \right] \right) 
\underset{n\to\infty}{\longrightarrow} z^{2}\psi'(1+\beta/2)/2
\] 
whose consequence is \eqref{dis}.  This, applied for $z=\pm 1$ together with Chebyshev inequality yields
\[
P(|n(E_{n}-\Delta_{n})-(\psi(1+\beta/2)-\log(\beta/2))|\ge\epsilon)\le Ce^{-\epsilon n^{1/2}}
\]
for a certain constant $C>0$.  This and an application of Borel-Cantelli's Lemma prove \eqref{con}.  
Again applying \eqref{k} with $n^{2}z$ in place of $z$, one can show that
\[
\frac{1}{n}\log\left(\E\left[ \exp\left( zn^{2}(E_{n}-\Delta_{n}) \right) \right] \right) 
\underset{n\to\infty}{\longrightarrow} R(z).
\]
for any $z\in\R$.  As a consequence of standard large deviations results (see for example Section 2.2 in \cite{S2}) we conclude the proof of the last part of the theorem. \qedhere
\end{proof}

\begin{corollary}
$E_{n}$ converges almost surely to $3/4$, the energy of the semicircular law on $[-2,2]$.  This implies that the spectral distribution, $\mu_{n}$ of $A_{n}$ converges almost surely to the semicircular law on $[-2,2]$.  
\end{corollary}

\begin{proof}
The convergence of $E_{n}$ to $3/4$ follows from \eqref{con} and the fact that the second expression in \eqref{e8} converges to $1/2-\int_{0}^{1}x\log(x)dx=3/4$.  Alternatively, we can use Proposition \ref{prop} for the convergence of $\Delta_{n}$ to the free entropy of the semicircular law.  For the converges of the spectral distribution, we use \ref{TQd} and \ref{prop} with $Q(x)=x^{2}/2$ plus the triangle inequality to justify that almost surely
\[
W_{2}(\mu_{n},\sigma)\le W_{2}(\mu_{n},\mu(\mathbf{y}_{n}))+W_{2}(\mu(\mathbf{y}_{n}),\sigma)\le \sqrt{2(E_{n}-\Delta_{n})}+\sqrt{2(3/4-\Delta_{n})}\underset{n\to\infty}{\longrightarrow}0.\qedhere
\]
\end{proof}

\bibliography{FreeTalagrand-corrected}
\bibliographystyle{mrl}

 \end{document}